# CONTACT PROCESSES ON RANDOM GRAPHS WITH POWER LAW DEGREE DISTRIBUTIONS HAVE CRITICAL VALUE 0

By Shirshendu Chatterjee and Rick Durrett[1]

*Cornell University*

If we consider the contact process with infection rate $\lambda$ on a random graph on $n$ vertices with power law degree distributions, mean field calculations suggest that the critical value $\lambda_c$ of the infection rate is positive if the power $\alpha > 3$. Physicists seem to regard this as an established fact, since the result has recently been generalized to bipartite graphs by Gómez-Gardeñes et al. [*Proc. Natl. Acad. Sci. USA* **105** (2008) 1399–1404]. Here, we show that the critical value $\lambda_c$ is zero for any value of $\alpha > 3$, and the contact process starting from all vertices infected, with a probability tending to 1 as $n \to \infty$, maintains a positive density of infected sites for time at least $\exp(n^{1-\delta})$ for any $\delta > 0$. Using the last result, together with the contact process duality, we can establish the existence of a quasi-stationary distribution in which a randomly chosen vertex is occupied with probability $\rho(\lambda)$. It is expected that $\rho(\lambda) \sim C\lambda^\beta$ as $\lambda \to 0$. Here we show that $\alpha - 1 \leq \beta \leq 2\alpha - 3$, and so $\beta > 2$ for $\alpha > 3$. Thus even though the graph is locally tree-like, $\beta$ does not take the mean field critical value $\beta = 1$.

**1. Introduction.** In this paper we will study the contact process on random graphs with a power-law degree distribution, i.e., for some constant $\alpha$, the degree of a typical vertex is $k$ with probability $p_k \sim Ck^{-\alpha}$ as $k \to \infty$. Following Newman, Strogatz and Watts (2001, 2002), we construct the random graph $G_n$ on the vertex set $\{1, 2, \ldots, n\}$ having degree distribution $\mathbf{p} = \{p_k : k \geq 0\}$ as follows. Let $d_1, \ldots, d_n$ be independent and have the distribution $P(d_i = k) = p_k$. We condition on the event $E_n = \{d_1 + \cdots + d_n \text{ is even}\}$ to have a valid degree sequence. As $P(E_n) \to 1/2$ as $n \to \infty$, the conditioning will have a little effect on the distribution of $d_i$'s. Having chosen the degree sequence $(d_1, d_2, \ldots, d_n)$, we allocate $d_i$ many half-edges to the vertex $i$, and

Received May 2008; revised February 2009.

[1]Supported in part by an NSF grant from the probability program.

then pair those half-edges at random. We also condition on the event that the graph is simple, i.e., it neither contains any self-loop at some vertex, nor contains multiple edges between two vertices. It can be shown [see e.g. Theorem 3.1.2 of Durrett (2007)] that if the degree distribution **p** has finite second moment, i.e., if $\alpha > 3$, the probability of the event that $G_n$ is simple has a positive limit as $n \to \infty$, and hence the conditioning on this event will not have much effect on the distribution of $d_i$'s.

We will be concerned with epidemics that take place on these random graphs. First consider the SIR (susceptible-infected-removed) model, in which sites begin as susceptible, and after being infected they get removed, i.e., become immune to further infection. In the simplest discrete-time formulation, an infected site $x$ at time $n$ will always be removed at time $n+1$ and for each susceptible neighbor $y$ at any time $n$, $x$ will cause $y$ to become infected at time $n+1$ with probability $p$, with all of the infection events being independent.

In this case the spreading of the epidemic is equivalent to percolation. To compute the threshold $p_c$ for a large, i.e., $O(n)$, epidemic to occur with positive probability, one notes that for a randomly chosen vertex $x$, the number of vertices at distance $m$ from $x$, $Z_m$, is approximately a two-phase branching process in which the number of first generation children has distribution **p**, but in the second and subsequent generations the offspring distribution is the size biased distribution $\mathbf{q} = \{q_k : k \geq 0\}$ satisfying

$$(1.1) \qquad q_{k-1} = \frac{kp_k}{\mu}, \qquad \text{where } \mu = \sum_k kp_k.$$

This occurs because vertices with degree $k$ are $k$ times as likely to be chosen for connections, and the edge that brings us to the new vertex uses up one of its degrees. For more details on this and the facts that we will quote in the next paragraph, see Chapter 3 of Durrett (2007).

With the above observation in hand, it is easy to compute the critical threshold for the SIR model. Let $\nu$ be the mean of the size biased distribution,

$$(1.2) \qquad \nu = \sum_k kq_k.$$

Suppose we start the infection at a randomly chosen vertex $x$. Now if $Y_m$ is the number of sites at distance $m$ from $x$ that become infected, then $EY_m = p\mu(p\nu)^{m-1}$. So the epidemic is supercritical if and only if $p > 1/\nu$. In particular, if $p_k \sim Ck^{-\alpha}$ as $k \to \infty$ and $\alpha \leq 3$, then $\nu = \infty$ and $p_c = 0$. Conversely if $\alpha > 3$ then $\nu < \infty$ and $p_c = 1/\nu > 0$. Hence for the SIR epidemic model on the random graph $G_n$ with power-law degree distribution, there is a positive threshold for the infection to survive if and only if the power $\alpha > 3$.

We will study the continuous-time SIS (susceptible-infected-susceptible) model and show that its behavior differs from that of the SIR model. In the SIS model, at any time $t$ each site $x$ is either infected or healthy (but susceptible). We often refer to the infected sites as occupied, and the healthy sites as vacant. We define the functions $\{\zeta_t : t \geq 0\}$ on the vertex set so that $\zeta_t(x)$ equals 0 or 1 depending on whether the site $x$ is healthy or infected at time $t$. An infected site becomes healthy at rate 1 independent of other sites and is again susceptible to the disease, while a susceptible site becomes infected at a rate $\lambda$ times the number of its infected neighbors. Harris (1974) introduced this model on the $d$-dimensional integer lattice and named it the *contact process.* See Liggett (1999) for an account of most of the known results. We will make extensive use of the *self-duality property* property of this process. If we let $\xi_t \equiv \{x : \zeta_t(x) = 1\}$ to be the set of infected sites at time $t$, we obtain a set-valued process. If we write $\xi_t^A$ to denote the process with $\xi_0^A = A$, then the self-duality property says that

(1.3) $$P(\xi_t^A \cap B \neq \varnothing) = P(\xi_t^B \cap A \neq \varnothing)$$

for any two subsets $A$ and $B$ of vertices.

Pastor-Satorras and Vespignani (2001a, 2001b, 2002) have made an extensive study of this model using mean-field methods. Their nonrigorous computations suggest the following conjectures about $\lambda_c$ the threshold for "prolonged persistence" of the contact process:

- If $\alpha \leq 3$, then $\lambda_c = 0$.
- If $3 < \alpha \leq 4$, then $\lambda_c > 0$ but the critical exponent $\beta$, which controls the rate at which the equilibrium density of infected sites goes to 0, satisfies $\beta > 1$.
- If $\alpha > 4$, then $\lambda_c > 0$ and the equilibrium density $\sim C(\lambda - \lambda_c)$ as $\lambda \downarrow \lambda_c$, i.e. the critical exponent $\beta = 1$.

Notice that the conjectured behavior of $\lambda_c$ for the SIS model parallels the results for $p_c$ in the SIR model quoted above.

Gómez-Gardeñes et al. (2008) have recently extended this calculation to the bipartite case, which they think of as a social network of sexual contacts between men and women. They define the polynomial decay rates for degrees in the two sexes to be $\gamma_M$ and $\gamma_F$, and argue that the epidemic is supercritical when the transmission rates for the two sexes satisfy

$$\sqrt{\lambda_M \lambda_F} > \lambda_c = \sqrt{\frac{\langle k \rangle_M \langle k \rangle_F}{\langle k^2 \rangle_F \langle k^2 \rangle_M}},$$

where the angle brackets indicate expected value and $k$ is shorthand for the degree distribution. Here $\lambda_c$ is positive when $\gamma_M, \gamma_F > 3$.

Our first goal is to show that $\lambda_c = 0$ for all $\alpha > 3$. Our proof starts with the following observation due to Berger et al. (2005). Here, we follow the formulation in Lemma 4.8.2 of Durrett (2007).

LEMMA 1.1. *Suppose $G$ is a star graph with center $0$ and leaves $1, 2, \ldots, k$. Let $A_t$ be the set of vertices infected in the contact process at time $t$ when $A_0 = \{0\}$. If $k\lambda^2 \to \infty$, then $P(A_{\exp(k\lambda^2/10)} \neq \varnothing) \to 1$.*

Based on results for the contact process on $(\mathbf{Z} \bmod n)$ by Durrett and Liu (1988) and Durrett and Schonmann (1988), and on $(\mathbf{Z} \bmod n)^d$ by Mountford (1993), it is natural to conjecture that in the contact process on $G_n$, with probability tending to 1 as $n \to \infty$, the infection survives for time $\geq \exp(cn)$ for some constant $c$. It certainly cannot last longer, because the total number of edges is $O(n)$, and so even if all sites are occupied at time 0, there is a constant $c$ so that with probability $\geq \exp(-cn)$ all sites will be vacant at time 1. Our next result falls a little short of that goal.

THEOREM 1. *Consider a Newman, Strogatz and Watts random graphs $G_n$ on the vertex set $\{1, 2, \ldots, n\}$, where the degrees $d_i$ satisfy $P(d_i = k) \sim Ck^{-\alpha}$ as $k \to \infty$ for some constant $C$ and some $\alpha > 3$, and $P(d_i \leq 2) = 0$. Let $\{\xi_t^1 : t \geq 0\}$ denote the contact process on the random graph $G_n$ starting from all sites occupied, i.e., $\xi_0^1 = \{1, 2, \ldots, n\}$. Then for any value of the infection rate $\lambda > 0$, there is a positive constant $p(\lambda)$ so that for any $\delta > 0$*

$$\inf_{t \leq \exp(n^{1-\delta})} P\left(\frac{|\xi_t^1|}{n} \geq p(\lambda)\right) \to 1 \qquad \text{as } n \to \infty.$$

One could assume that $\nu > 1$ and look at the process on the giant component, but we would rather avoid this complication. The assumption $P(d_i \leq 2) = 0$ is convenient, because it implies the following.

LEMMA 1.2. *Consider a Newman, Strogatz and Watts graphs, $G_n$, on $n$ vertices, where the degrees of the vertices, $d_i$, satisfy $P(d_i \leq 2) = 0$, and the mean of the size biased degree distribution $\nu < \infty$. Then*

$$P(G_n \text{ is connected}) \to 1 \qquad \text{as } n \to \infty,$$

*and if $D_n$ is the diameter of $G_n$,*

$$P(D_n > (1 + \varepsilon) \log n / \log \nu) \to 0 \qquad \text{for any } \varepsilon > 0.$$

The size of the giant component in the graph is given by the nonextinction probability of the two-phase branching process, so $P(d_i \leq 2) = 0$ is needed to have the size $\sim n$. Intuitively, Lemma 1.2 is obvious because the worst case is the random 3-regular graph, and in this case, the graph is not only connected and has diameter $\sim (\log n)/(\log 2)$, see Sections 7.6 and 10.3 of Bollobás (2001), but the probability of a Hamiltonian cycle tends to 1, see Section 9.3 of Janson, Luczak, and Ruciński (2000). We have not been able to find a proof of Lemma 1.2 in the literature, so we give one in Section 5. By comparing the growth of the cluster with a branching process it is easy to show $P(D_n < (1-\varepsilon)\log n/\log \nu) \to 0$ for any $\varepsilon > 0$.

In a sense the main consequence of Theorem 1 is not new. Berger et al. (2005), see also (2009), show that $\lambda_c = 0$ for a generalization of the Bárabasi–Albert model in which each new point has $m$ edges which are with probability $\beta$ connected to a vertex chosen uniformly at random and with probability $1-\beta$ to a vertex chosen with probability proportional to its degree. Theorem 2 in Cooper and Frieze (2003) shows that such graphs have power law degree distributions with $\alpha = 1+2/(1-\beta)$, so these examples have $\alpha \in [3, \infty)$ and $\lambda_c = 0$.

Having acknowledged the previous work of BBCS, it should be noted that (i) our result applies to a large class of power law graphs that have a different structure; and (ii) the BBCS proof yields a lower bound on the persistence time of $\exp(cn^{1/(\alpha-1)})$ compared to our $\exp(n^{1-\delta})$. Our improved bound on the survival times relies only on the power law degree distribution and the fact that the diameter is bounded by $C\log n$, so it also applies to graphs BBCS consider.

Theorem 1 shows that the fraction of infected sites in the graph $G_n$ is bounded away from zero for a time longer than $\exp(n^{1/2})$. So using self-duality we can now define a quasi-stationary measure $\xi_\infty^1$ on the subsets of $\{1, 2, \ldots, n\}$ as follows. For any subset of vertices $A$, $P(\xi_\infty^1 \cap A \neq \varnothing) \equiv P(\xi_{\exp(n^{1/2})}^A \neq \varnothing)$. Let $X_n$ be uniformly distributed on $\{1, 2, \ldots, n\}$ and let $\rho_n(\lambda) = P(X_n \in \xi_\infty^1)$. Berger et al. (2005) show that for the contact process on their preferential attachment graphs, there are positive, finite constants so that

$$b\lambda^C \leq \rho_n(\lambda) \leq B\lambda^c.$$

In contrast, we get reasonably good numerical bounds on the critical exponent.

THEOREM 2. *Suppose $\alpha > 3$. There is a $\lambda_0 > 0$ so that if $0 < \lambda < \lambda_0$ and $0 < \delta < 1$, then there exists two constants $c(\alpha, \delta)$ and $C(\alpha, \delta)$ so that as $n \to \infty$*

$$P(c\lambda^{1+(\alpha-2)(2+\delta)} \leq \rho_n(\lambda) \leq C\lambda^{1+(\alpha-2)(1-\delta)}) \to 1.$$

When $\alpha$ is close to 3 and $\delta$ is small, the powers in the lower and upper bounds are close to 3 and 2. The ratio of the two powers is $\leq (2+\delta)/(1-\delta) \approx 2$ when $\delta$ is small.

The intuition behind the lower bound is that if the infection starts from a vertex of degree $d(x) \geq (10/\lambda)^{2+\delta}$, then it survives for a long time with a probability bounded away from 0. The density of such points is $C\lambda^{(2+\delta)(\alpha-1)}$, but we can improve the bound to the one given by looking at neighbors of these vertices, which have density $C\lambda^{(2+\delta)(\alpha-2)}$ and will infect their large degree neighbor with probability $\geq c\lambda$.

For the upper bound we show that if $m(\alpha,\delta)$ is large enough and the infection starts from a vertex $x$ such that there is no vertex of degree $\geq \lambda^{-(1-\delta)}$ within distance $m$ from $x$, then its survival is very unlikely. To get the extra factor of $\lambda$ we note that the first event must be a birth. Based on the proof of Lemma 1.1, we expect that survival is unlikely if there is no nearby vertex of degree $\geq \lambda^{-2}$ and hence the lower bound gives the critical exponent.

It is natural to speculate that the density of the quasi-stationary measure $\rho_n(\lambda) \to \rho(\lambda)$ as $n \to \infty$. By the heuristics for the computation of $\lambda_c$ in the SIR model, it is natural to guess that, when $\alpha > 2$, $\rho(\lambda)$ is the expected probability of weak survival for the contact process on a tree generated by the two-phase branching process, starting with the origin occupied.

Here the phrase 'weak survival' refers to set of infected sites being not empty for all times, in contrast to 'strong survival' where the origin is re-infected infinitely often. As in the case of the contact process on the Bollobás–Chung small world studied by Durrett and Jung (2007), it is the weak survival critical value that is the threshold for prolonged persistence on the finite graph.

*Sketch of the proof of Theorem 1.* The remainder of the paper is devoted to proofs. Let $V_n^\varepsilon$ be the set of vertices in the graph $G_n$ with degree at least $n^\varepsilon$. We call the points in $V_n^\varepsilon$ *stars*. We say that a star of degree $k$ is *hot* if at least $\lambda k/4$ of its neighbors are infected and is *lit* if at least $\lambda k/10$ of its neighbors are infected. Our first step, taken in Lemma 2.2, is to improve the proof of Lemma 1.1 to show that a hot star will remain lit for time $\exp(cn^\varepsilon)$ with high probability.

To keep the system going for a long time, we cannot rely on just one star. There are $O(n^{1-\varepsilon(\alpha-1)})$ stars in this graph which has diameter $O(\log n)$. If one star goes out, presence of a lit star can make it hot again within a time $2n^{\varepsilon/3}$ with probability at least $n^{-b}$. See Lemmas 2.3 and 2.4 for this. Lemma 2.6 shows that a lit star gets hot within $2\exp(n^{\varepsilon/3})$ units of time with probability

$$\geq 1 - 5\exp(-\lambda^2 n^{\varepsilon/3}/16),$$

and Lemma 2.5 shows that a hot star eventually succeeds to make a nonlit star hot within $\exp(n^{\varepsilon/2})$ units of time with probability

$$\geq 1 - 8e^{-\lambda^2 n^\varepsilon/80}.$$

Using these estimates, we can show that the number of lit stars dominates a random walk with a strong positive drift, and hence more than 3/4's of the collection will stay lit for a time $O(\exp(n^{1-\alpha\varepsilon}))$. See Proposition 1 at the end of Section 2 for the argument.

To get a lower bound on the density of infected sites, first we bound the probability of the event that the dual process, starting from a vertex of degree $(10/\lambda)^{2+\delta}$, reaches more than 3/4's of the stars. We do this in two steps. In the first step (see Lemma 3.2) we get a lower bound for the probability of the dual process reaching one of the stars. To do this, we consider a chain of events in which we reach vertices with degree $(10/\lambda)^{k+\delta}$ for $k \geq 2$ sequentially. In the second step (see Lemma 3.3) we again use a comparison with random walk to show that, with probability tending to 1, the dual process, starting from any lit star, will light up more than 3/4's of the stars. Then we show that the above events are asymptotically uncorrelated, and use a second moment argument to complete the proof of Theorem 1 and the lower bound for the density in Theorem 2.

OPEN PROBLEM. Improve the bounds in Theorem 2 and extend the result to $\alpha > 1$.

When $2 < \alpha < 3$ the size biased distribution has infinite mean. Chung and Lu (2002, 2003) obtained bounds on the diameter in this case, and later van der Hofstadt, Hooghiemstra and Zamenski (2007) showed if $H_n$ is the distance between 1 and 2 then

$$H_n \sim \frac{2 \log \log n}{-\log(\alpha - 2)}.$$

When $1 < \alpha < 2$ the size-biased distribution has infinite mass. van der Esker et al. (2005) have shown in this case

$$\lim_{n \to \infty} P(H_n = 2) = \lim_{n \to \infty} 1 - P(H_n = 3) = p \in (0, 1)$$

so the graph is very small.

All of the results about the persistence of infection at stars in Section 2 are valid for any $\alpha$, since they only rely on properties of the contact process on a star graph and an upper bound on the diameter. The results in Section 3, rely on the existence of the size biased distribution and hence are restricted to $\alpha > 2$. The proof of the lower bound should be extendible to that case, but the proof of the upper bound given in Section 4 relies heavily on the size-biased distribution having finite mean. When $1 < \alpha < 2$, the size-biased distribution does not exist and the situation changes drastically. We guess that in this case $\rho_n(\lambda) = O(\lambda)$.

**2. Persistence of infection at stars.** Let $\varepsilon > 0$ and let $V_n^\varepsilon$ be the set of vertices in our graph $G_n$ with degree at least $n^\varepsilon$. We call these vertices *stars*. We say that a vertex of degree $k$ is *hot* if it has at least $L = \lambda k/4$ infected neighbors and we call it *lit* if it has at least $0.4L = \lambda k/10$ infected neighbors. We will show that if $\varepsilon$ is small, then in the contact process starting from all vertices occupied, most of the stars in $V_n^\varepsilon$ will remain lit for time $O(\exp(n^{1-\alpha\varepsilon}))$.

We begin with a slight improvement of Lemma 1.1 which gives a numerical estimate of the failure probability, but before that we need two simple estimates.

LEMMA 2.1. *If $0 \le x \le a \le 1$ then $e^x \le 1 + (1+a)x$ and $e^{-x} \le 1 - (1 - 2a/3)x$.*

PROOF. Using the series expansion for $e^x$

$$e^x \le 1 + x + \frac{ax}{2}\left(1 + \frac{1}{2} + \left(\frac{1}{2}\right)^2 + \cdots\right),$$

$$e^{-x} \le 1 - x + \frac{ax}{2}\left(1 + \left(\frac{1}{2}\right)^2 + \left(\frac{1}{2}\right)^4 + \cdots\right)$$

and summing the geometric series gives the result. □

LEMMA 2.2. *Let $G$ be a star graph with center $0$ and leaves $1, 2, \ldots, k$. Let $A_t$ be the set of vertices infected in the contact process at time $t$. Suppose $\lambda \le 1$ and $\lambda^2 k \ge 50$. Let $L = \lambda k/4$ and let $T = \exp(k\lambda^2/80)/4L$. Let $P_{L,i}$ denote the probability when at time $0$ the center is at state $i$ and $L$ leaves are infected. Then*

$$P_{L,i}\left(\inf_{t \le T} |A_t| \le 0.4L\right) \le 7e^{-\lambda^2 k/80} \qquad \text{for } i = 0, 1.$$

PROOF. Write the state of the system as $(m, n)$ where $m$ is the number of infected leaves and $n = 1$ if the center is infected and $0$ otherwise. To reduce to a one-dimensional chain, we will concentrate on the first coordinate. When the state is $(m, 0)$ with $m > 0$, the next event will occur after exponential time with mean $1/(m\lambda + m)$, and the probability that it will be the reinfection of the center is $\lambda/(\lambda + 1)$. So the number of leaf infections $N$ that will die while the center is $0$ has a shifted geometric distribution with success probability $\lambda/(\lambda + 1)$, i.e.,

$$P(N = j) = \left(\frac{1}{\lambda + 1}\right)^j \cdot \frac{\lambda}{\lambda + 1} \qquad \text{for } j \ge 0.$$

Let $N_L$ be the realization of $N$ when the state of the system is $(L, 0)$. Then $N_L$ will be more than $0.1L$ with probability

$$(2.1) \quad P_{L,0}(N_L > 0.1L) \leq (1+\lambda)^{-0.1L} \leq e^{-\lambda L/20} = e^{-\lambda^2 k/80}.$$

Here we use the inequality $1 + \lambda \geq e^{\lambda/2}$. If $N_L \leq 0.1L$, then there will be at least $0.9L$ infected leaves when the center is infected.

The next step is to modify the chain so that the infection rate is 0 when the number of infected leaves is $L = \lambda k/4$ or greater. In this case the number of infected leaves $\geq Y_t$ where

$$\begin{aligned} & & \text{at rate} \\ Y_t &\to Y_t - 1, & \lambda k/4, \\ Y_t &\to Y_t + 1, & 3\lambda k/4 \text{ for } Y_t < L, \\ Y_t &\to Y_t - N, & 1. \end{aligned}$$

To bound the survival time of this chain, we will estimate the probability that starting from $0.8L$ it will return to $0.4L$ before hitting $L$. During this time $Y_t$ is a random walk that jumps at rate $\lambda k + 1$. Let $X$ be the change in the random walk in one step. Then

$$X = \begin{cases} -1, & \text{with probability } (\lambda k/4)/(\lambda k + 1), \\ +1, & \text{with probability } (3\lambda k/4)/(\lambda k + 1), \\ -N, & \text{with probability } 1/(\lambda k + 1), \end{cases}$$

and so

$$Ee^{\theta X} = e^\theta \cdot \frac{3}{4} \cdot \frac{\lambda k}{\lambda k + 1} + e^{-\theta} \cdot \frac{1}{4} \cdot \frac{\lambda k}{\lambda k + 1}$$
$$+ \frac{1}{\lambda k + 1} \sum_{j=0}^\infty e^{-\theta j} \left(\frac{1}{\lambda + 1}\right)^j \cdot \frac{\lambda}{\lambda + 1}.$$

If $e^{-\theta}/(\lambda + 1) < 1$, the third term on the right is

$$\frac{\lambda}{\lambda k + 1} \cdot \frac{1}{1 + \lambda - e^{-\theta}}.$$

If we pick $\theta < 0$ so that $e^{-\theta} = 1 + \lambda/2$, then

$$Ee^{\theta X} = \frac{\lambda k}{\lambda k + 1} \left(\frac{1}{1 + \lambda/2} \cdot \frac{3}{4} + (1 + \lambda/2) \cdot \frac{1}{4} + \frac{2}{\lambda k}\right).$$

Since $1/(1 + x) < 1 - x + x^2$ for $0 < x < 1$,

$$\frac{1}{1 + \lambda/2} \cdot \frac{3}{4} + (1 + \lambda/2) \cdot \frac{1}{4} + \frac{2}{\lambda k} - 1 < \left(-\frac{\lambda}{2} + \frac{\lambda^2}{4}\right)\frac{3}{4} + \frac{\lambda}{8} + \frac{2}{\lambda k}$$
$$< -\frac{3\lambda}{16} + \frac{\lambda}{8} + \frac{2}{\lambda k},$$

where in the last inequality, we have used $\lambda < 1$. Since we have assumed $\lambda^2 k \geq 50$, the right-hand side is $< 0$.

To estimate the hitting probability we note that if $\phi(x) = \exp(\theta x)$ and $Y_0 \geq 0.6L$, then $\phi(Y_t)$ is a supermartingale until it hits $L$. Let $q$ be the probability that $Y_t$ hits the interval $(-\infty, 0.4L]$ before returning to $L$. Since $\theta < 0$, we have $\phi(x) \geq \phi(0.4L)$ for $x \leq 0.4L$. So using the optional stopping theorem we have

$$q\phi(0.4L) + (1-q)\phi(L) \leq \phi(0.8L),$$

which implies that

$$q \leq \phi(0.8L)/\phi(0.4L) = \exp(0.4\theta L) \leq e^{-\lambda^2 k/40},$$

as $e^{-\theta} = 1 + \lambda/2 \geq e^{\lambda/4}$ when $\lambda/4 < 1/2$ (sum the series for $e^x$).

At this point we have estimated the probability that the chain started at a point $\geq 0.8L$ will go to $L$ before going below $0.4L$. When the chain is at $L$, the time until the next jump is exponential with mean $1/(L+1) \geq 1/2L$. The probability that the jump takes us below $0.8L$ is (since $1 + \lambda \geq e^{\lambda/2}$)

$$\leq (1+\lambda)^{-0.2L} \leq e^{-\lambda L/10} = e^{-\lambda^2 k/40}.$$

Thus the probability that the chain fails to return to $L$, $M = e^{\lambda^2 k/80}$ times before going below $0.4L$ is

$$\leq 2e^{-\lambda^2 k/80}.$$

Using Chebyshev's inequality on the sum, $S_M$ of $M$ exponentials with mean 1 (and hence variance 1),

$$P(S_M < M/2) \leq 4/M.$$

Multiplying by $1/2L$ we see that the total time, $T_M$ of the first $M$ excursions satisfies

$$P(T_M < M/4L) \leq 4e^{-\lambda^2 k/80}.$$

Combining this with the previous estimate on the probability of having fewer than $M$ returns and the error probability in (2.1) proves the desired result. □

Thus Lemma 2.2 shows that a hot star will remain lit for a long time with probability very close to 1. Our next step is to investigate the process of transferring the infection from one star to another. The first step in doing that is to estimate what happens when only the center of the star infected.

LEMMA 2.3. *Let $G$ be a star graph with center 0 and leaves $1, 2, \ldots, k$. Let $0 < \lambda < 1$, $\delta > 0$ and suppose $\lambda^{2+\delta} k \geq 10$. Again let $P_{l,i}$ denote the probability when at time 0 the center is in state $i$ and $l$ leaves are infected. Let $\tau_0$ be the first time 0 becomes healthy, and let $T_j$ be the first time the number of infected leaves equals $j$. If $L = \lambda k/4$, $\gamma = \delta/(4+2\delta)$, and $K = \lambda k^{1-\gamma}/4$, then for $k \geq k_0(\delta)$*

$$P_{0,1}(T_K > \tau_0) \leq 2/k^\gamma,$$
$$P_{K,1}(T_0 < T_L) \leq \exp(-\lambda^2 k^{1-\gamma}/16) \leq 1/k^\gamma,$$
$$E_{0,1}(T_L \mid T_L < \infty) \leq 2.$$

Combining the first two inequalities $P_{0,1}(T_L < \infty) \geq 1 - 2/k^\gamma$, and using Markov's inequality, if we can infect a vertex of degree at least $k$ such that $k \geq k_0(\delta)$ and $\lambda^{2+\delta} k > 10$, then with probability $\geq 1 - 5/k^\gamma$ the vertex gets hot within the next $k^\gamma$ units of time.

PROOF OF LEMMA 2.3. Note that $\tau_0 \sim \exp(1)$, and for any $t \leq \tau_0$, the leaves independently becomes healthy at rate 1 and infected at rate $\lambda$. Let $p_0(t)$ is the probability that leaf $j$ is infected at time $t$ when the central vertex of the star has remained infected for all times $s \leq t$. $p_0(0) = 0$ and

$$\frac{dp_0(t)}{dt} = -p_0(t) + (1 - p_0(t))\lambda = \lambda - (\lambda + 1)p_0(t).$$

So solving gives $p_0(t) = \int_0^t \lambda e^{-(\lambda+1)(t-s)}\, ds = \frac{\lambda}{\lambda+1}(1 - e^{-(\lambda+1)t})$. From this it follows that

(2.2) $\quad P_{0,1}(T_K < \tau_0) \geq P(\text{Binomial}(k, p_0(k^{-\gamma})) > K) P(\tau_0 > k^{-\gamma}).$

Now if $k^\gamma > 8/3$, $(\lambda+1)k^{-\gamma} \leq 3/4$ and it follows from Lemma 2.1 that

$$p_0(k^{-\gamma}) \geq \lambda k^{-\gamma}/2.$$

Writing $p = p_0(k^{-\gamma})$ to simplify formulas, if $\theta > 0$

$$P(\text{Binomial}(k, p) \leq K) \leq e^{\theta K}(1 - p + pe^{-\theta})^k.$$

Since $\log(1+x) \leq x$ the right-hand side is

$$\leq \exp\left(\frac{\theta \lambda k^{1-\gamma}}{4} + (e^{-\theta} - 1)\frac{\lambda k^{1-\gamma}}{2}\right).$$

Taking $\theta = 1/2$ and using Lemma 2.1 to conclude $e^{-1/2} - 1 \leq -1/3$, the above is

$$\leq \exp(-\lambda k^{1-\gamma}/24) \leq \exp(-k^{1/2-\gamma}/8),$$

since $\lambda^2 k \geq 9$. Using this in (2.2), the right-hand side is

$$\geq (1 - \exp(-k^{1/2-\gamma}/8))(1 - k^{-\gamma}) \geq 1 - 2/k^\gamma,$$

if $k^{1/2-\gamma} \geq 8\gamma \log k$.

Using the supermartingale from the proof of Lemma 2.2, if $q = P_{K,1}(T_0 < T_L)$, then we have

$$q \cdot 1 + (1-q)e^{\theta L} \leq e^{\theta K},$$

and so $q \leq e^{\theta K} \leq e^{-\lambda K/4}$. In the last step we have used $e^\theta = 1/(1 + \lambda/2) \leq e^{-\lambda/4}$, which comes from Lemma 2.1. Filling in the value of $K$, $e^{-\lambda K/4} = e^{-\lambda^2 k^{1-\gamma}/16}$. Now

$$\lambda^2 k^{1-\gamma} = (\lambda^{2+\delta} k)^{2/(2+\delta)} k^{1-\gamma-2/(2+\delta)} \geq 10^{2/(2+\delta)} k^{\delta/(4+2\delta)}.$$

So if $k^{\delta/(4+2\delta)} > 16 \cdot 10^{-2/(2+\delta)} \gamma \log k$, then $e^{-\lambda K/4} \leq 1/k^\gamma$.

To bound the time we use the lower bound random walk $Y_t$ from Lemma 2.2. $EN = 1/\lambda$, so

$$EY_t = \left(\frac{\lambda k}{2} - \frac{1}{\lambda}\right)t = \left(\frac{\lambda^2 k - 2}{2\lambda}\right)t.$$

Let $T_L^Y$ be the hitting time of $L$ for the random walk $Y_t$. Using the optional stopping theorem for the martingale $Y_t - (\lambda^2 k - 2)t/2\lambda$ and the bounded stopping time $T_L^Y \wedge t$ we get

$$EY_{T_L^Y \wedge t} - \left(\frac{\lambda^2 k - 2}{2\lambda}\right)E(T_L^Y \wedge t) = EY_0 = 0.$$

Since $EY_{T_L^Y \wedge t} \leq L = \lambda k/4$, it follows that

$$E(T_L^Y \wedge t) \leq \left(\frac{2\lambda}{\lambda^2 k - 2}\right)L = \frac{\lambda^2 k/2}{\lambda^2 k - 2} = \frac{1}{2 - 4/\lambda^2 k} \leq 1,$$

as by our assumption $\lambda^2 k \geq 4$. Letting $t \to \infty$ we have $ET_L^Y \leq 1$. Since $Y_t$ is a lower bound for the number of infected leaves, $T_L \mathbf{1}_{[T_L < \infty]} \leq T_L^Y$. Hence

$$E_{0,1}(T_L | T_L < \infty) = \frac{E_{0,1}(T_L \mathbf{1}_{[T_L < \infty]})}{P_{0,1}(T_L < \infty)}$$

$$\leq \frac{E_{0,1} T_L^Y}{P_{0,1}(T_K < \tau_0) P_{K,1}(T_L < T_0)} \leq \frac{1}{1/2} = 2$$

for large $k$. $\square$

To transfer infection from one vertex to another we use the following lemma.

LEMMA 2.4. *Let $v_0, v_1, \ldots, v_m$ be a path in the graph and suppose that $v_0$ is infected at time 0. Then the probability that $v_m$ will become infected by time $m$ is $\geq (e^{-1}(1-e^{-\lambda})e^{-1})^m$.*

PROOF. The first factor is the probability that the infection at $v_0$ lasts for time 1, the second the probability that $v_0$ infects $v_1$ by time 1, and the third the probability that the infection at $v_1$ remains until time 1. Iterating this $m$ times brings the infection from 0 to $m$. □

When the diameter of the graph is $\leq 2\log n$, the probability in Lemma 2.4 is $\geq n^{-b}$ for some $b \in (1/2, \infty)$, and the time required is $\leq 2\log n$. Combining this with Lemma 2.3 (with $k = n^\varepsilon$ and $\gamma = 1/3$) shows that if $n$ is large, then with probability $\geq Cn^{-b}$ we can use one hot star to make another star hot within time $2n^{\varepsilon/3}$. Using Lemma 2.2 and trying repeatedly gives the following lemma.

LEMMA 2.5. *Let $s_1$ and $s_2$ be two stars in $V_n^\varepsilon$ and suppose that $s_1$ is hot at time 0. Then, for large $n$, $s_2$ will be hot by time $T = \exp(n^{\varepsilon/2})$ with probability*

$$\geq 1 - 8e^{-\lambda^2 n^\varepsilon/80}.$$

PROOF. If $n$ is large, Lemma 2.2 shows that $s_1$ remains lit for $T$ units of time with probability $\geq 1 - 7e^{-\lambda^2 n^\varepsilon/80}$. Let $t_n = 2n^{\varepsilon/3}$ and consider the discrete time points $t_n, 2t_n, \ldots$. At all of these time points we can think of a path starting from an infected neighbor of $s_1$ up to $s_2$. Using one such path the infection gets transmitted to $s_2$ and it gets hot in $2n^{\varepsilon/3}$ units of time with probability $\geq Cn^{-b}$ for some constant $C$. So $s_1$ fails to make $s_2$ hot by time $T$ with probability

$$\leq (1 - Cn^{-b})^{T/t_n} \leq \exp(-Cn^{-b}T/t_n) \leq \exp(-\lambda^2 n^\varepsilon/80)$$

for large $n$. For the first inequality we use $1 - x \leq e^{-x}$. Combining with the first error probability in this proof, we get the result. □

Next we show that a lit star becomes hot with a high probability, and then helps to make other nonlit stars lit.

LEMMA 2.6. *Let $s$ be a star of $V_n^\varepsilon$ and suppose that $s$ is lit at time 0. Then $s$ will be hot by time $2\exp(n^{\varepsilon/3})$ with probability*

$$\geq 1 - 5\exp(-\lambda^2 n^{\varepsilon/3}/16) \qquad \text{if } n \text{ is large.}$$

PROOF. Since $s$ is lit, it has at least $\lambda n^{\varepsilon}/10$ infected neighbors at time 0. If $s$ itself is not infected at time 0, let $N$ be the number of leaf infections that die out before $s$ gets infected. Using similar argument as in the beginning of the proof of Lemma 2.2,

$$P(N = j) = \left(\frac{1}{\lambda + 1}\right)^j \cdot \frac{\lambda}{\lambda + 1} \qquad \text{for } j \geq 0,$$

which implies

$$P(N > \lambda n^{\varepsilon}/20) \leq (1+\lambda)^{-\lambda n^{\varepsilon}/20} \leq e^{-\lambda^2 n^{\varepsilon}/40},$$

as $1 + \lambda > e^{\lambda/2}$ by Lemma 2.1. Also the time $T_M$ taken for $M = \lambda n^{\varepsilon}/20$ leaf infections to die out is a sum of $M$ exponentials with mean at most $1/(\lambda+1)M \leq 1/M$. Now if $n^{2\varepsilon/3} > 40/16$, the above error probability is $\leq e^{-\lambda^2 n^{\varepsilon/3}/16}$.

Using Chebyshev's inequality on the sum, $S_M$ of $M$ exponentials with mean 1 (and hence variance 1), we see that if $\exp(n^{\varepsilon/3}) \geq 2$, i.e., $n^{\varepsilon/3} > \log 2$

$$P(S_M > M \exp(n^{\varepsilon/3})) \leq \frac{1}{M(\exp(n^{\varepsilon/3}) - 1)^2}$$

$$\leq \frac{4}{M \exp(2n^{\varepsilon/3})} \leq \exp(-\lambda^2 n^{\varepsilon/3}/16),$$

where in the final inequality we have used $M > 4$, i.e., $n^{\varepsilon} > 80/\lambda$, and $\lambda^2/16 < 2$.

Multiplying by $1/M$ we see that the total time, $T_M$, satisfies

$$P(T_M > \exp(n^{\varepsilon/3})) \leq \exp(-\lambda^2 n^{\varepsilon/3}/16).$$

Combining these two error probabilities gives that $s$ will be infected along with at least $\lambda n^{\varepsilon}/20$ infected neighbors within $\exp(n^{\varepsilon/3})$ units of time with error probability

(2.3) $$\leq 2 \exp(-\lambda^2 n^{\varepsilon/3}/16).$$

Now $\lambda n^{\varepsilon}/20 \geq \lambda n^{\varepsilon/3}/4$, when $n^{2\varepsilon/3} > 5$. So if $s$ is infected and has at least $\lambda n^{\varepsilon}/20$ infected neighbors, then using the second inequality of Lemma 2.3 (with $\gamma = 2/3$ and $k = n^{\varepsilon}$), $s$ becomes hot with error probability

$$\leq \exp(-\lambda^2 n^{\varepsilon/3}/16).$$

Finally using Markov's inequality and the third inequality of Lemma 2.3, the time $T_s$ taken by $s$ to get hot, after it became infected, is more than $T = \exp(n^{\varepsilon/3})$ with probability

$$\leq 2\exp(-n^{\varepsilon/3}) \leq 2\exp(-\lambda^2 n^{\varepsilon/3}/16),$$

as $\lambda < 1$. Combining all these error probabilities proves the lemma. □

We now use Lemmas 2.5, 2.6 and 2.2 to prove that if the contact process starts from all sites infected, then for a long time at least 3/4's of the stars will be lit.

PROPOSITION 1. *Let $I_{n,t}^\varepsilon$ be the set of stars in $V_n^\varepsilon$ which are lit at time $t$ in the contact process $\{\xi_t^1 : t \geq 0\}$ on $G_n$. Let $t_n = 2\exp(n^{\varepsilon/2})$ and $M_n = \exp(n^{1-\alpha\varepsilon})$. Then there is a stopping time $T_n$ such that $T_n \geq M_n \cdot t_n$ and*

$$P(|I_{n,T_n}^\varepsilon| \leq (3/4)|V_n^\varepsilon|) \leq \exp(-Cn^\varepsilon).$$

PROOF. Let $\alpha_n = |V_n^\varepsilon|$. Clearly $|I_{n,0}^\varepsilon| = \alpha_n$. We will estimate the probability that starting from $(7/8)\alpha_n$ lit stars, the number goes below $(3/4)\alpha_n$ before reaching $\alpha_n$. Define the stopping times $\tau_i$s' and $\sigma_i$s' as follows. Let $\tau_0 = \sigma_0 = 0$ and for $i \geq 0$ let

$$\tau_{i+1} \equiv \inf\{t > \tau_i + \sigma_i t_n : |I_{n,t}^\varepsilon| = (7/8)\alpha_n\},$$
$$\sigma_{i+1} \equiv \min\{s \in \mathbb{N} : |I_{n,\tau_{i+1}+s\cdot t_n}^\varepsilon| \notin ((3/4)\alpha_n, \alpha_n)\}.$$

We need to look at time lags that are multiples of $t_n$ in the definition of $\sigma_i$ because in our worst nightmare (which is undoubtedly a paranoid delusion) all the lit stars of degree $k \geq n^\varepsilon$ at time $\tau_{i+1}$ have exactly $0.1k$ infected neighbors.

Lemma 2.6 implies that a lit star of $V_n^\varepsilon$ gets hot within time $2\exp(n^{\varepsilon/3}) \leq \exp(n^{\varepsilon/2})$ (for large $n$) with probability $\geq 1 - 5\exp(-\lambda^2 n^{\varepsilon/3}/16)$. Combining with Lemma 2.2 gives that a lit star at time 0 gets hot by time $t_n/2$ and remains lit at time $t_n$ with probability $\geq 1 - 6\exp(-\lambda^2 n^{\varepsilon/3}/16)$ for large $n$. Now if $|I_{n,t}^\varepsilon| < \alpha_n$ for some $t$, then the number of lit stars will increase at time $t + t_n$ with probability $\geq P(A \cap B)$, where

- A: All the lit stars will get hot by $t_n/2$ units of time, and be lit for time $t_n$.
- B: A nonlit star will become hot by time $t_n/2$ in presence of another hot star, and remain lit for another $t_n/2$ units of time.

Now using the above argument $P(A) \geq 1 - 6n\exp(-\lambda^2 n^{\varepsilon/3}/16)$, as there are at most $n$ stars. Combining Lemmas 2.5 and 2.2 gives $P(B) \geq 1 - 9\exp(-\lambda^2 n^\varepsilon/80)$. So $P(A \cap B) \geq 1 - \exp(-n^{\varepsilon/4})$ for large $n$. Using the stopping times $|I_{n,\tau_i+r\cdot t_n}^\varepsilon| \geq W_r$ for $r \leq \sigma_i$, where $\{W_r : r \geq 0\}$ is a discrete time random walk satisfying

(2.4)
$$W_r \to W_r - 1 \quad \text{with probability } \exp(-n^{\varepsilon/4}),$$
$$W_r \to W_r + 1 \quad \text{with probability } 1 - \exp(-n^{\varepsilon/4}),$$

and $W_0 = (7/8)\alpha_n$. Now $\theta^{W_r}$ is a martingale where

(2.5)
$$\theta = \frac{\exp(-n^{\varepsilon/4})}{1 - \exp(-n^{\varepsilon/4})} < \exp(-n^{\varepsilon/4}/2).$$

If $q$ is the probability that $W_r$ goes below $(3/4)\alpha_n$ before hitting $\alpha_n$, then applying the optional stopping theorem

$$q \cdot \theta^{(3/4)\alpha_n} + (1-q) \cdot \theta^{\alpha_n} \leq \theta^{(7/8)\alpha_n},$$

which implies

$$q \leq \theta^{(\alpha_n/8)} \leq \exp(-Cn^{1-(\alpha-1)\varepsilon}),$$

as $\alpha_n \sim Cn^{1-(\alpha-1)\varepsilon}$ for some constant $C$. So the probability that the random walk fails to return to $\alpha_n$ at least $M_n = \exp(n^{1-\alpha\varepsilon})$ times before going below $(3/4)\alpha_n$ is $\leq \exp(-Cn^\varepsilon)$. Now if

$$K = \min\{i \geq 1 : |I^\varepsilon_{n,\tau_i+\sigma_i \cdot t_n}| \leq (3/4)\alpha_n\},$$

the coupling with the random walk will imply $P(K \leq M_n) \leq \exp(-Cn^\varepsilon)$, and hence for $T_n \equiv \tau_{M_n} + \sigma_{M_n} \cdot t_n$

$$P(|I^\varepsilon_{n,T_n}| \leq (3/4)|V^\varepsilon_n|) \leq \exp(-Cn^\varepsilon).$$

As $\sigma_i \geq 1$ for all $i$, by our construction $T_n \geq M_n \cdot t_n$, and we get the result. □

So the infection persists for time longer than $\exp(n^{1-\alpha\varepsilon})$ in the stars of $V^\varepsilon_n$.

**3. Density of infected stars.** Proposition 1 implies that if the contact process starts with all vertices infected, most of the stars remain lit even after $\exp(n^{1-\alpha\varepsilon})$ units of time. In this section we will show that the density of infected stars is bounded away from 0 and we will find a lower bound for the density. We start with the following lemma about the growth of clusters in the random graph $G_n$, when we expose the neighbors of a vertex one at a time. For more details on this procedure see Section 3.2 of Durrett (2007).

LEMMA 3.1. *Suppose $0 < \delta \leq 1/8$. Let $A$ be the event that the two clusters, starting from 1 and 2 respectively, intersect before their sizes grow to $n^\delta$. Then*

$$P(A) \leq Cn^{-(1/4-\delta)}.$$

PROOF. If $d_1, \ldots, d_n$ are the degrees of the vertices, then

$$(3.1) \qquad P\Big(\max_{1 \leq i \leq n} d_i > n^{3/(2\alpha-2)}\Big) \leq n \cdot P(d_1 > n^{3/(2\alpha-2)}) \leq c/\sqrt{n}$$

for some constant $c$. Suppose all the degrees are at most $n^{3/(2\alpha-2)}$. Suppose $R_1$ and $R_2$ are the clusters starting from 1 and 2 up to size $n^\delta$. Let $B$ be the event that $R_1$ contains a vertex of degree $\geq n^{1/(2\alpha-2)}$. Let $e_n$ be the sum of degrees of all those vertices with degree $\geq n^{1/(2\alpha-2)}$. While growing $R_1$ the probability that a vertex of degree $\geq n^{1/(2\alpha-2)}$ will be included on any step is

$$\leq \frac{e_n}{\sum_{i=1}^{n} d_i - n^{\delta+3/(2\alpha-2)}} \equiv \beta_n.$$

Since the size biased distribution is $q_k \sim Ck^{-(\alpha-1)}$ as $k \to \infty$, $\sum_{s \geq k} q_s \sim Ck^{-(\alpha-2)}$ as $k \to \infty$, and we have $e_n \sim Cn^{1-(\alpha-2)/(2\alpha-2)}$ and hence $\beta_n \sim Cn^{-(\alpha-2)/(2\alpha-2)}$ as $n \to \infty$. So for large $n$ $\beta_n \leq c_1 n^{-1/4}$ for some constant $c_1$, when $\alpha > 3$. Thus

$$P(B^c) \geq 1 - c_1/n^{1/4-\delta}.$$

If $B^c$ occurs, all the degrees of the vertices of $R_1$ are at most $n^{1/(2\alpha-2)}$. In that case, while growing $R_2$ the probability of choosing one vertex from $R_1$ is

$$\leq \frac{n^{\delta+1/(2\alpha-2)}}{\sum_{i=1}^{n} d_i - n^{\delta+3/(2\alpha-2)}} \leq c_2/n^{1-\delta-1/(2\alpha-2)}.$$

So the conditional probability

$$P(A^c|B^c) \geq (1 - c_2 n^{-(1-\delta-1/(2\alpha-2))})^{n^\delta} \geq 1 - c_2/n^{1-2\delta-1/(2\alpha-2)}.$$

Hence combining these two

$$P(A^c) \geq (1 - c_1/n^{1/4-\delta})(1 - c_2/n^{1-2\delta-1/(2\alpha-2)}) \geq 1 - C_1/n^{1/4-\delta},$$

and that completes the proof. □

Lemma 3.1 will help us to show that in the dual contact process, staring from any vertex of degree $\geq (10/\lambda)^{2+\delta}$ for some $\delta > 0$, the infection reaches a star of $V_n^\varepsilon$, with probability bounded away from 0.

LEMMA 3.2. *Let $\xi_t^A$ be the contact process on $G_n$ starting from $\xi_0^A = A$. Suppose $0 < \varepsilon < 1/20(\alpha - 1)$. Then there are constants $\lambda_0 > 0$, $n_0 < \infty$, $c_0 = c_0(\lambda, \varepsilon)$ and $p_i > 0$ independent of $\lambda < \lambda_0$, $n \geq n_0$ and $\varepsilon$ such that if $T = n^{c_0}$, $v_2$ is a vertex with degree $d(v_2) \geq (10/\lambda)^{2+\delta}$ for some $0 < \delta < 1$ and $v_1$ is a neighbor of $v_2$,*

$$P(\xi_T^{\{v_2\}} \cap V_n^\varepsilon) \geq p_2, \qquad P(\xi_{T+1}^{\{v_1\}} \cap V_n^\varepsilon) \geq p_1 \lambda.$$

PROOF. The second conclusion follows immediately from the first, since the probability that $v_1$ will infect $v_2$ before time 1, and that $v_2$ will stay infected until time 1 is

$$\geq \frac{\lambda}{\lambda+1}(1 - e^{-(\lambda+1)})e^{-1} \geq c\lambda.$$

Let $\Lambda_m$ be the set of vertices in $G_n$ of degree $\geq (10/\lambda)^{m+\delta}$ for $m \geq 2$. Define $\gamma = \frac{\delta}{2(2+\delta)}$ and

$$B = 2(\alpha-1)\log(10/\lambda), \qquad u = (e^{-1}(1-e^{-\lambda})e^{-1})^{-(B+1)},$$

$$w_n \equiv \log(n^\varepsilon)/\log(10/\lambda) - \delta, \qquad T_m = T_m^1 + T_m^2,$$

where $T_m^1 = (10/\lambda)^{(m+\delta)\gamma}$, $T_m^2 = u^m$, and we let $n^{c_0} = \sum_{m=2}^{w_n} T_m$.

Define $E_2 = \{\xi_{T_2}^{\{v_2\}} \cap \Lambda_3 \neq \varnothing\}$ and for $m \geq 3$, having defined $E_2, \ldots, E_{m-1}$, we let

$$E_m = \{\xi_{T_m}^{\{v_m\}} \cap \Lambda_{m+1} \neq \varnothing\} \quad \text{and} \quad v_m \in \xi_{T_{m-1}}^{\{v_{m-1}\}} \cap \Lambda_m.$$

Let $A_m$ be the event that the clusters of size $(10/\lambda)^{(m+\delta+1)(\alpha-2)}$ starting from two neighbors of $v_m$ do not intersect and

$$F = \bigcap_{m=2}^{w_n} A_m.$$

Since $\varepsilon < 1/20(\alpha - 1)$, the cluster size $(10/\lambda)^{(m+\delta+1)(\alpha-2)}$ is at most $n^{1/10}$ for $m \leq w_n$. So using Lemma 3.1 and the fact $\binom{k}{2} < k^2$,

$$P(F^c) \leq \left( \sum_{m=2}^{w_n} (10/\lambda)^{2m+2\delta} \right) cn^{-(1/4-1/10)}$$

$$\leq n^{2\varepsilon} cn^{-(1/4-1/10)} < cn^{-(1/4-3/20)} < 1/6$$

for large $n$.

Since each vertex has degree at least 3, if $F$ occurs then by the choice of $B$ the neighborhood of radius $Bm$ around $v_m$ will contain more than $(\frac{10}{\lambda})^{(m+\delta+1)(\alpha-2)+m}$ vertices. Let $G_m$ be the event that the neighborhood of radius $Bm$ around $v_m$ intersects $\Lambda_{m+1}$. In the neighborhood of $v_m$ probability of having a vertex of $\Lambda_{m+1}$ is at least $c(\lambda/10)^{(m+\delta+1)(\alpha-2)}$. Hence

$$P(G_m^c \cap F) \leq (1 - c(\lambda/10)^{(m+\delta+1)(\alpha-2)})^{(10/\lambda)^{m+(m+\delta+1)(\alpha-2)}}$$

$$\leq \exp(-(10/\lambda)^m).$$

If $\lambda$ is small, $\sum_{m=2}^{\infty} \exp(-(10/\lambda)^m) \leq 1/6$.

On the intersection of $F$ and $G_m$ we have a vertex of $\Lambda_{m+1}$ within radius $Bm$ of $v_m$. Using Lemmas 2.2 and 2.3, in the contact process $\{\xi_t^{\{v_m\}} : t \geq 0\}$, $v_m$ gets hot at time $T_m^1$ and remains lit till time $T_m$ with error probability $\leq c\lambda^{(m+\delta)\gamma}$ for small $\lambda$. If $v_m$ is lit, then Lemma 2.4 shows that $v_m$ fails to transfer the infection to some vertex in $\Lambda_{m+1}$ within time $T_m^2$ with probability

$$\leq [1 - (e^{-1}(1-e^{-\lambda})e^{-1})^{Bm}]^{T_m^2/(Bm)}$$

$$\leq \exp[-(e^{-1}(1-e^{-\lambda})e^{-1})^{-m}/(Bm)] \equiv \eta_m,$$

where $\equiv$ indicates we are making a definition, and hence $P(E_m^c G_m F) \leq c\lambda^{(m+\delta)\gamma} + \eta_m$. If $\lambda$ is small $\sum_{m=2}^{w_m}[c\lambda^{(m+\delta)\gamma} + \eta_m] \leq 1/6$, we can take $p_2 = 1/2$ and the proof is complete. $\square$

Lemma 3.2 gives a lower bound on the probability that an infection starting from a neighbor of a vertex of degree $\geq (10/\lambda)^{2+\delta}$ reaches a star. Lemma 2.3 shows that if the infection reaches a star, then with probability tending to 1 the star gets hot within $n^{\varepsilon/3}$ units of time. Combining these two we get the following.

PROPOSITION 2. *Suppose $0 < \varepsilon < 1/20(\alpha-1)$. There are constants $\lambda_0 > 0$, $n_0 < \infty$ $c_1 = c_1(\lambda, \varepsilon)$ and $p_1 > 0$, which does not depend on $\lambda < \lambda_0$, $n \geq n_0$ and $\varepsilon$, such that for any vertex $v_1$ with a neighbor $v_2$ of degree $d(v_2) \geq (10/\lambda)^{2+\delta}$ for some $\delta \in (0,1)$, and $T = n^{c_1}$ the probability that $\xi_T^{\{v_1\}}$ contains a hot star is bounded below by $p_1\lambda$.*

Next we will show that if we start with one lit star, then after time $\exp(n^{\varepsilon/2})$ at least 3/4's of the stars will be lit.

LEMMA 3.3. *Let $I_{n,t}^\varepsilon$ be the set of stars which are lit at time $t$ in the contact process on $G_n$ such that $|I_{n,0}^\varepsilon| = 1$. Then for $T' = \exp(n^{\varepsilon/2})$*

$$P(|I_{n,T'}^\varepsilon| < (3/4)|V_n^\varepsilon|) \leq 7\exp(-\lambda^2 n^{\varepsilon/3}/16).$$

PROOF. Let $s_1$ be the lit star at time 0. As seen in Proposition 1, $s_1$ remains lit at time $T' = \exp(n^{\varepsilon/2})$ with probability $\geq 1 - 6\exp(-\lambda^2 n^{\varepsilon/3}/16)$ for large $n$. With probability $\geq Cn^{-b}$ another star gets hot within time $t_n = 2n^{\varepsilon/3}$ and remains lit at time $T'$. Using similar argument as in Lemma 2.5, the process fails to make $(3/4)|V_n^\varepsilon|$ many stars lit by time $T'$ with probability

$$\leq (3/4)|V_n^\varepsilon|(1 - Cn^{-b})^{T'/t_n}$$
$$\leq (3/4)|V_n^\varepsilon|\exp(-Cn^{-b}T'/t_n) \leq \exp(-\lambda^2 n^{\varepsilon/3}/16),$$

as $|V_n^\varepsilon| = Cn^{1-(\alpha-1)\varepsilon}$ and $1 - x \leq e^{-x}$. So combining with the earlier error probability we get the result. $\square$

Now we are almost ready to prove our main result. However, we need one more lemma that we will use in the proof of the theorem.

LEMMA 3.4. *Let $F$ and $G$ be two events which involve exposing $n^\delta$ many vertices starting at 1 and 2 respectively for some $0 < \delta \leq 1/8$. Then*

$$|P(F \cap G) - P(F)P(G)| \leq Cn^{-(1/4-\delta)}.$$

PROOF. Let $R_1$ and $R_2$ be the clusters for exposing $n^\delta$ many vertices starting from 1 and 2 respectively, and let $A$ be the event that they intersect. Clearly

$$P(F \cap G) \leq P(A) + P(F \cap G \cap A^c)$$
$$= P(A) + P(F \cap A^c)P(G \cap A^c)$$
$$\leq P(A) + P(F)P(G).$$

Using similar argument for $F^c$ and $G$ we get

$$|P(F \cap G) - P(F)P(G)| \leq P(A).$$

We estimate $P(A)$ using Lemma 3.1. □

Lemma 3.4 shows that two events which involve exposing at most $n^{1/8}$ vertices starting from two different vertices are asymptotically uncorrelated. Now we give the proof of the main theorem.

*Proof of Theorem 1.* Given $\delta > 0$, choose $\varepsilon = \min\{\delta/\alpha, 1/20(\alpha - 1)\}$. Let $A_n$ be the set of vertices in $G_n$ with a neighbor of degree at least $(10/\lambda)^{2+\delta}$. Clearly $|A_n|/n \to c_0(\lambda/10)^{(2+\delta)(\alpha-2)}$ as $n \to \infty$ for some constant $c_0$. Define the random variables $Y_x, x \in A_n$ as $Y_x = 1$ if the dual contact process starting from $x$ can light up a star of $V_n^\varepsilon$ and 0 otherwise. By Proposition 2, $EY_x \geq p_1\lambda$ for some constant $p_1 > 0$ and for any $x \in A_n$.

If we grow the cluster starting from $x \in A_n$ and exposing one vertex at a time, we can find a star on any step with probability at least $cn^{-(\alpha-2)\varepsilon}$. So with probability $1 - \exp(-cn^\varepsilon)$, we can find a star of $V_n^\varepsilon$ within the exposure of at most $n^{\alpha\varepsilon}$ vertices. So, with high probability, lighting a star up is an event involving at most $n^{(\alpha+1)\varepsilon}$ many vertices. As $(\alpha+1)\varepsilon < 1/8$, using Lemma 3.4, we can say

$$P(Y_x = 1, Y_z = 1) - P(Y_x = 1)P(Y_z = 1)$$
$$\leq (1 - \exp(-cn^\varepsilon))Cn^{-(1/4-(\alpha+1)\varepsilon)} + \exp(-cn^\varepsilon) \equiv \theta_n.$$

Using our bound on the covariances,

$$\text{var}\left(\sum_{x \in A_n} Y_x\right) \leq n + \binom{n}{2}\theta_n,$$

and Chebyshev's inequality gives

$$P\left(\left|\sum_{x \in A_n} (Y_x - EY_x)\right| \geq n\gamma\right) \leq \frac{n + \binom{n}{2}\theta_n}{n^2\gamma^2} \to 0 \quad \text{as } n \to \infty$$

for any $\gamma > 0$, since $\theta_n \to 0$ as $n \to \infty$. Since $EY_x \geq p_1 \lambda$ and $|A_n|/n \to c_0(\lambda/10)^{(2+\delta)(\alpha-2)}$, if we take $p_l \equiv p_1 \lambda \cdot c_0(\lambda/10)^{(2+\delta)(\alpha-2)}/2$ then

$$(3.2) \quad \lim_{n \to \infty} P\bigg(\sum_{x \in A_n} Y_x \geq np_l\bigg) = 1.$$

Now if $Y_x = 1$, Proposition 2 says that the dual process starting from $x$ makes a star hot after $T_1 = n^{c_1}$ units of time. Then by Lemma 3.3 within next $T_2 = \exp(n^{\varepsilon/2})$ units of time the dual process lights up 75% of all the stars with probability $1 - 7\exp(-\lambda^2 n^{\varepsilon/3}/16)$.

Let $I_{n,t}^{\varepsilon}$ be the set of stars which are lit at time $t$ in the contact process $\{\xi_t^1 : t \geq 0\}$ and

$$T_3 = \inf\{t > \exp(n^{1-\alpha\varepsilon}) : |I_{n,t}^{\varepsilon}| \geq (3/4)|V_n^{\varepsilon}|\}.$$

By Proposition 1, $P(T_3 < \infty) \geq 1 - \exp(-cn^{\varepsilon})$. Let

$$\mathcal{S} = \{S \subset \{1, 2, \ldots, n\} : \xi_t^1 = S \Rightarrow |I_{n,t}^{\varepsilon}| \geq (3/4)|V_n^{\varepsilon}|\}.$$

Using the Markov property and self-duality of the contact process we get the following inequality. For any subset $B$ of the vertex set, and for the event $F_n \equiv [T_3 < \infty]$ we have

$$P[(\xi_{T_1+T_2+T_3}^{\mathbf{1}} \supset B) \cap F_n]$$
$$= \sum_{S \in \mathcal{S}} P(\xi_{T_1+T_2}^S \supset B) P(\xi_{T_3}^1 = S | F_n) P(F_n)$$
$$= \sum_{S \in \mathcal{S}} P(\xi_{T_1+T_2}^{\{x\}} \cap S \neq \varnothing \ \forall x \in B) P(\xi_{T_3}^1 = S | F_n) P(F_n)$$
$$\geq \sum_{S \in \mathcal{S}} P(|\xi_{T_1+T_2}^{\{x\}} \cap I_{n,T_3}^{\varepsilon}| > (3/4)|V_n^{\varepsilon}| \ \forall x \in B) P(\xi_{T_3}^1 = S | F_n) P(F_n)$$
$$\geq P(Y_x = 1 \ \forall x \in B)(1 - 7|B|\exp(-\lambda^2 n^{\varepsilon/3}/16))P(F_n)$$
$$\geq P(Y_x = 1 \ \forall x \in B)(1 - 2\exp(-cn^{\varepsilon/4})),$$

as $|B| \leq n$ and $P(F_n) \geq 1 - \exp(-cn^{\varepsilon})$. Hence for $T = T_1 + T_2 + T_3$, combining with (3.2) and using the attractiveness property of the contact process we conclude that as $n \to \infty$

$$(3.3) \quad \inf_{t \leq T} P\bigg(\frac{|\xi_t^{\mathbf{1}}|}{n} > p_l\bigg) = P\bigg(\frac{|\xi_T^{\mathbf{1}}|}{n} > p_l\bigg)$$
$$\geq P\bigg(\xi_T^{\mathbf{1}} \supseteq \{x : Y_x = 1\}, \sum_{x \in A_n} Y_x \geq np_l\bigg) \to 1,$$

which completes the proof of Theorem 1, and proves the lower bound in Theorem 2.

**4. Upper bound in Theorem 2.** For the upper bound, we will show that if the infection starts from a vertex $x$ with no vertex of degree $> 1/\lambda^{1-\delta}$ nearby, it has a very small chance to survive. To get the 1 in upper bound we need to use the fact that first event in the contact process starting at $x$ has to be a birth so we begin with that calculation.

Let $\Lambda_\delta$ be the set of vertices of degree $> \lambda^{\delta-1}$. Define $Z_x, x \in \{1, 2, \ldots, n\}$ as $Z_x = 1$ if the dual contact process $\{\xi_t^{\{x\}} : t \geq 0\}$ starting from $x$ survives for $T' = 1/\lambda^{\alpha-1}$ units of time, and 0 otherwise. We will show $EZ_x \leq C\lambda^{1+(\alpha-2)(1-\delta)}$ for some constant $C$. If $T_1$ is the time for the first event in the dual process, then $ET_1 \leq 1$ and using Markov's inequality $P(T_1 > 1/\lambda^{\alpha-1}) < \lambda^{\alpha-1}$. So if $T_1 < 1/\lambda^{\alpha-1}$, the first event must be a birth for $Z_x$ to be 1. So for $x \in \Lambda_\delta$,

$$P(Z_x = 1) \leq P(T_1 > 1/\lambda^{\alpha-1}) + \sum_{i > \lambda^{\delta-1}} p_i \frac{\lambda i}{\lambda i + 1}$$

$$\leq \lambda^{\alpha-1} + C\lambda \sum_{i > \lambda^{\delta-1}} i^{-(\alpha-1)}$$

$$\leq \lambda^{\alpha-1} + C\lambda \cdot \lambda^{(\alpha-2)(1-\delta)}.$$

For $x \in \Lambda_\delta^c$, let $w(\lambda) \leq C\lambda^{(\alpha-2)(1-\delta)}$ be the size-biased probability of having a vertex of $\Lambda_\delta$ in its neighborhood. If $d(x) = i$, the expected number of vertices in a radius $m$ around $x$ is at most $i \cdot EZ_m$, where $Z_m$ is the total progeny up to $m$th generation of the branching process with offspring distribution $q_k = (k+1)p_{k+1}/\mu \sim ck^{\alpha-1}$. So the expected number of vertices, which are within a distance $m = \lceil(\alpha-1)/\delta\rceil$, the smallest integer larger than $(\alpha-1)/\delta$, from $x$ and belong to $\Lambda_\delta$, is

$$\leq \sum_{i=2}^{(1/\lambda)^{1-\delta}} p_i \cdot i \cdot EZ_m \cdot C\lambda^{(\alpha-2)(1-\delta)} \leq C\lambda^{(\alpha-2)(1-\delta)}.$$

Using Markov's inequality the probability of having at least one vertex of $\Lambda_\delta$ within a distance $m$ from $x$ has the same upper bound as above.

Until we reach $\Lambda_\delta$, $|\xi_t^{\{x\}}| \leq Y_t$ where

$$Y_t \to Y_t - 1 \quad \text{at rate } Y_t,$$
$$Y_t \to Y_t + 1 \quad \text{at rate } Y_t \lambda \cdot (1/\lambda)^{1-\delta} = Y_t \lambda^\delta.$$

So $Y_t$ jumps at rate $Y_t(1 + \lambda^\delta)$ and it jumps to $Y_t + 1$ with probability $\lambda^\delta/(1 + \lambda^\delta) < \lambda^\delta$. If $T_1 < 1/\lambda^{\alpha-1}$, the first event in the dual process $\xi_t^{\{x\}}$ must be a birth for $Z_x$ to be 1. Let $T_{2m}$ is the time of the $2m$th event after the first event. Then $ET_{2m} \leq 2m/(1 + \lambda^\delta)$ and using Markov's inequality

$$P(T_{2m} > 1/\lambda^{\alpha-1}) \leq C\lambda^{\alpha-1}.$$

Now if $T_{2m} < 1/\lambda^{\alpha-1}$ and there is no vertex of $\Lambda_\delta$ within a distance $m$ of $x$, the infection starting at $x$ survives for time $T'$ only if $Y_t$ has at least $m$ up jumps before hitting 0. If there are $\leq m-1$ up jumps in the first $2m$ then $Y_t$ will hit 0 by $T_{2m}$, as $Y_0 = 2$. The probability of this event is

$$\leq P(B \geq m), \quad \text{where } B \sim \text{Binomial}(2m, \lambda^\delta),$$
$$\leq 2^{2m}\lambda^{m\delta} \leq 2^{2m}\lambda^{\alpha-1}.$$

Combining all three error probabilities, for any $x \in \Lambda_\delta^c$,

$$P(Z_x = 1) \leq P(T_1 > 1/\lambda^{\alpha-1}) + P(T_{2m} > 1/\lambda^{\alpha-1})$$
$$+ \sum_{i \leq \lambda^{\delta-1}} p_i \frac{\lambda i}{\lambda i + 1} \cdot C\lambda^{(\alpha-2)(1-\delta)}$$
$$\leq C\lambda^{1+(\alpha-2)(1-\delta)}.$$

Using an argument similar to one at the end of the proof of Theorem 1

$$P\left(\left|\sum_x (Z_x - EZ_x)\right| > n\gamma\right) \to 0 \quad \text{as } n \to \infty$$

for any $\gamma > 0$. Since $EZ_x \leq C\lambda^{1+(\alpha-2)(1-\delta)}$ for all $x \in \{1, 2, \ldots, n\}$, if we take $p_u = 3C\lambda^{1+(\alpha-2)(1-\delta)}$, then

$$P\left(\sum_x Z_x \geq np_u\right) \to 0 \quad \text{as } n \to \infty.$$

So by making $C$ larger in the definition of $p_u$ and using the attractiveness of the contact process

$$\inf_{t \geq T'} P(|\xi_t^1| \leq p_u n) \to 1$$

as $n \to \infty$.

**5. Proof of connectivity and diameter.** We conclude the paper with the proof of Lemma 1.2. We begin with a large deviations result. The fact is well-known, but the proof is short so we give it for completeness.

LEMMA 5.1. *Let $X_1, X_2, \ldots$ be i.i.d., nonnegative with mean $\mu$. If $\rho < \mu$, then there is a constant $\gamma > 0$ so that*

$$P(X_1 + \cdots + X_k \leq \rho k) \leq e^{-\gamma k}.$$

PROOF. Let $\phi(\theta) = Ee^{-\theta X}$. If $\theta > 0$ then

$$e^{-\theta \rho k} P(X_1 + \cdots + X_k \leq \rho k) \leq \phi(\theta)^k.$$

So we have

$$P(X_1 + \cdots + X_k \leq \rho k) \leq \exp(k\{\theta\rho + \log \phi(\theta)\}).$$

$\log(\phi(0)) = 0$ and as $\theta \to 0$

$$\frac{d}{d\theta} \log(\phi(\theta)) = \frac{\phi'(\theta)}{\phi(\theta)} \to -\mu.$$

So $\log \phi(\theta) \sim -\mu\theta$ as $\theta \to 0$, and the result follows by taking $\theta$ small. □

PROOF OF LEMMA 1.2. We will prove the result in the following steps:

*Step* 1. Let $k_n = (\log n)^2$. The size of the cluster $C_x$, starting from $x \in \{1, 2, \ldots, n\}$, reaches size $k_n$ with probability $1 - o(n^{-1})$.

*Step* 2. There is a $B < \infty$ so that if the size of $C_x$ reaches size $B \log n$, it will reach $n^{2/3}$ with probability $1 - O(n^{-2})$.

*Step* 3. Let $\zeta > 0$. Two clusters $C_x$ and $C_y$, starting from $x$ and $y$ respectively, of size $n^{(1/2)+\zeta}$ will intersect with probability $1 - o(n^{-2})$.

Steps 2 and 3 follow from the proof of Theorem 3.2.2 of Durrett (2007), so it is enough to do Step 1. Before doing this, note that if $d_1, \ldots, d_n$ are the degrees of the vertices, and $\eta > 0$ then as $n \to \infty$,

$$P\Big(\max_{1 \leq i \leq n} d_i > n^{(1+\eta)/(\alpha-1)}\Big) \leq n \cdot P(d_1 > n^{(1+\eta)/(\alpha-1)}) \sim C n^{-\eta}.$$

Given $\alpha > 3$, we choose $\eta > 0$ small enough so that $(1+\eta)/(\alpha-1) < 1/2$.

To prove Step 1, we will expose one vertex at a time. Following the notation of Durrett (2007), suppose $A_t, U_t$ and $R_t$ are the sets of active, unexplored and removed sites respectively at time $t$ in the process of growing the cluster starting from 1, with $R_0 = \{1\}$, $A_0 = \{z : 1 \sim z\}$ and $U_0 = \{1, 2, \ldots, n\} - A_0 \cup R_0$. At time $\tau = \inf\{t : A_t = \varnothing\}$ the process stops. If $A_t \neq \varnothing$, pick $i_t$ from $A_t$ in some way measurable with respect to the process up to that time and let

$$R_{t+1} = R_t \cup \{i_t\},$$
$$A_{t+1} = A_t \cup \{z \in U_t : i_t \sim z\} - \{i_t\},$$
$$U_{t+1} = U_t - \{z \in U_t : i_t \sim z\}.$$

Here $|R_t| = t+1$ for $t \leq \tau$ and so $C_1 = \tau + 1$. If there were no collisions, then $|A_{t+1}| = |A_t| - 1 + Z$ where $Z$ has the size biased degree distribution $q$. Let $q^\eta$ be the distribution of $(Z \mid Z \leq n^{(1+\eta)/(\alpha-1)})$. Then on the event $\{\max_i d_i \leq n^{(1+\eta)/(\alpha-1)}\}$, $|A_t|$ is dominated by a random walk $S_t = S_0 + Z_1 + \cdots + Z_t$, where $S_0 = A_0$ and $Z_i \sim q^\eta$. Since $q_{k-1} = kp_k/\mu$, we have $q_0 = q_1 = 0$ and hence $q_0^\eta = q_1^\eta = 0$. Then $S_t$ increases monotonically.

If we let $T = \inf\{m : S_m \geq k_n\}$ then

(5.1) $\quad P(|C_1| \leq k_n) \leq P(S_t - |A_t| \geq 4 \text{ for some } t \leq T).$

As observed above, if $n$ is large, all of the vertices have degree $\leq n^\beta$ where $\beta = (1+\eta)/(\alpha-1) < 1/2$. As long as $S_t \leq 2k_n$, each time we add a new vertex and the probability that it is in the active set is at most

$$\gamma_n = \frac{2k_n n^\beta}{\sum_{i=1}^n d_i - 2k_n n^\beta} \leq C k_n n^{\beta-1}$$

for large $n$. Thus the probability of two or more collisions while $S_t \leq 2k_n$ is $\leq (2k_n)^2 \gamma_n^2 = o(n^{-1})$.

If $S_T - S_{T-1} \leq k_n$, then the previous argument suffices, but $S_T - S_{T-1}$ might be as large as $n^\beta$. Letting $m > 1/(1-2\beta)$, we see that the probability of $m$ or more collisions is at most

$$(n^\beta)^m (C n^{\beta-1})^m = o(n^{-1}).$$

To grow the cluster we will use a breadth first search: we will expose all the vertices at distance 1 from the starting point, then those at distance 2, etc. When a collision occurs, we do not add a vertex, and we delete the one with which a collision has occurred, so two are lost. There is at most one collision while $S_t \leq 2k_n$. Since $S_0 \geq 3$, it is easy to see that the worst thing that can happen in terms of the growth of the cluster is for the collision to occur on the first step, reducing $S_0$ to 1. After this the number of vertices doubles at each step so size $k_n$ is reached before we have gone a distance $\log_2 k_n$ from the starting point.

In the final step we might have a jump $S_\tau - S_{\tau-1} \geq k_n$ and $m$ collisions, but as long as $k_n = (\log n)^2 > 2m$ we do not lose any ground. In the growth before time $T$, each vertex, except for possibly one collision, has added two new vertices to the active set. From this it is easy to see that the number of vertices in the active set is at least $k_n/2 - 2m$.

To grow the graph now, we will expose all of the vertices in the current active set, then expose all of the neighbors of these vertices, etc. Let $\varepsilon > 0$. The proof of Theorem 3.2.2 in Durrett (2007) shows (see page 78) that if $\delta$ is small then until $n\delta$ vertices have been exposed, the cluster growth dominates a random walk with mean $\nu - \varepsilon$. Let $J_1, J_2, \ldots$ be the successive sizes of the active set when these phases are complete. The large deviations result, Lemma 5.1, implies that there is a $\gamma > 0$ so that

$$P(J_{i+1} \leq (\nu - 2\varepsilon) J_i \mid J_i = j_i) \leq \exp(-\gamma j_i).$$

Since $J_1 \geq (\log n)^2/2 - 8$, it follows from this result that with probability $\geq 1 - o(n^{-1})$, in at most

$$\left(\frac{1}{2} + \zeta\right) \frac{\log n}{\log(\nu - \varepsilon)}$$

steps, the active set will grow to size $n^{(1/2)+\zeta}$. Using the result from Step 3 and noting that the initial phase of the growth has diameter $\leq \log_2 k_n = O(\log \log n)$ the desired result follows. $\square$

## REFERENCES

Berger, N., Borgs, C., Chayes, J. T. and Saberi, A. (2005). On the spread of viruses on the internet. In *Proceedings of the Sixteenth Annual ACM-SIAM Symposium on Discrete Algorithms* 301–310 (electronic). ACM, New York. MR2298278

Berger, N., Borgs, C., Chayes, J. T. and Saberi, A. (2009). Weak local limits for preferential attachment graphs. To appear.

Bollobás, B. (2001). *Random Graphs*, 2nd ed. *Cambridge Studies in Advanced Mathematics* **73**. Cambridge Univ. Press, Cambridge. MR1864966

Chung, F. and Lu, L. (2002). The average distances in random graphs with given expected degrees. *Proc. Natl. Acad. Sci. USA* **99** 15879–15882 (electronic). MR1944974

Chung, F. and Lu, L. (2003). The average distance in a random graph with given expected degrees. *Internet Math.* **1** 91–113. MR2076728

Cooper, C. and Frieze, A. (2003). A general model of web graphs. *Random Structures Algorithms* **22** 311–335. MR1966545

Durrett, R. (2007). *Random Graph Dynamics*. Cambridge Univ. Press, Cambridge. MR2271734

Durrett, R. and Jung, P. (2007). Two phase transitions for the contact process on small worlds. *Stochastic Process. Appl.* **117** 1910–1927. MR2437735

Durrett, R. and Liu, X. F. (1988). The contact process on a finite set. *Ann. Probab.* **16** 1158–1173. MR942760

Durrett, R. and Schonmann, R. H. (1988). The contact process on a finite set. II. *Ann. Probab.* **16** 1570–1583. MR958203

Gómez-Gardeñes, J., Latora, V., Moreno, Y. and Profumo, E. (2008). Spreading of sexually transmitted diseases in heterosexual populations. *Proc. Natl. Acad. Sci.* **105** 1399–1404.

Harris, T. E. (1974). Contact interactions on a lattice. *Ann. Probab.* **2** 969–988. MR0356292

Janson, S., Łuczak, T. and Rucinski, A. (2000). *Random Graphs*. Wiley, New York. MR1782847

Liggett, T. M. (1999). *Stochastic Interacting Systems: Contact, Voter and Exclusion Processes. Grundlehren der Mathematischen Wissenschaften [Fundamental Principles of Mathematical Sciences]* **324**. Springer, Berlin. MR1717346

Mountford, T. S. (1993). A metastable result for the finite multidimensional contact process. *Canad. Math. Bull.* **36** 216–226. MR1222537

Newman, M. E. J., Strogatz, S. H. and Watts, D. J. (2001). Random graphs with arbitrary degree distributions and their applications. *Phys. Rev. E.* **64** 026118.

Newman, M. E. J., Strogatz, S. H. and Watts, D. J. (2002). Random graph models of social networks. *Proc. Nat. Acad. Sci.* **99** 2566–2572.

Pastor-Satorras, R. and Vespignani, A. (2001a). Epidemic spreading in scale-free networks. *Phys. Rev. Letters* **86** 3200–3203.

Pastor-Satorras, R. and Vespignani, A. (2001b). Epidemic dynamics and endemic states in complex networks. *Phys. Rev. E.* **63** 066117.

Pastor-Satorras, R. and Vespignani, A. (2002). Epidemic dynamics in finite size scale-free networks. *Phys. Rev. E.* **65** 035108(R).

van den Esker, H., van der Hofstad, R., Hooghiemstra, G. and Znamenski, D. (2005). Distances in random graphs with infinite mean degrees. *Extremes* **8** 111–141 (2006). MR2275914

van der Hofstad, R., Hooghiemstra, G. and Znamenski, D. (2007). Distances in random graphs with finite mean and infinite variance degrees. *Electron. J. Probab.* **12** 703–766 (electronic). MR2318408

Operations Research
  and Information Engineering
Rhodes Hall
Cornell University
Ithaca, New York 14853
USA
E-mail: sc499@cornell.edu

Mathematics
Malott Hall
Cornell University
Ithaca, New York 14853
USA
E-mail: rtd1@cornell.edu